\documentclass[12pt]{article}
\usepackage{amsmath,amssymb}
\usepackage{rotating}
\usepackage{xypic}
\usepackage{amssymb}
\usepackage{amsmath,amscd}
\usepackage{pst-node,pstricks,multido,pst-plot,pst-text,pst-3d}%
\usepackage{graphicx}
\usepackage{amsfonts}
\newtheorem{example}{Example}[section]

\newtheorem{remark}[example]{Remark}
\newtheorem{theorem}[example]{Theorem}

\newtheorem{corollary}[example]{Corollary}

\newtheorem{proposition}[example]{Proposition}

\newtheorem{lemma}[example]{Lemma}

\def\Proof{\noindent \it Proof -- \rm}
\def\qed{\hspace{3.5mm} \hfill \vbox{\hrule height 3pt depth 2 pt width 2mm}
\bigskip\\}

\def\S{{\mathfrak  S}}
\def\cal#1{{\mathfrak #1}}

\def\<{\langle}
\def\>{\rangle}

\def\C{{\mathbb C}}

\def\ie{{\it i.e.}, }

\def\goth{\mathfrak}

\def\ashuff#1#2#3{
\kern 1pt \vrule height#1 \overline{\vrule height#3 width 0pt
\hskip#2} \rule{.3pt}{#1}\overline{\vrule height#3 width 0pt
\hskip#2} \rule{.3pt}{#1} \kern 1pt }

\def\Det{{\rm Det}}
\def\Id{{\rm Id}}

\def\DetF{{\rm Det}_{\goth F}}

\def\sign{\mbox{sign}}

\begin{document}

\title{Hyperdeterminants on semilattices}
\author{Jean-Gabriel Luque
\thanks{\centerline{\tiny
Institut Gaspard Monge,
 Universit\'e de Marne-la-Vall\'ee,}\newline\centerline{\tiny
 77454 Marne-la-Vall\'ee cedex 2}\newline\centerline{\tiny
  France}\newline\centerline{\tiny
email:Jean-Gabriel.Luque@univ-mlv.fr}}}

\markboth{Hyperdeterminants on semilattices}{Hyperdeterminants on
semilattices}
 \maketitle
\begin{abstract}
 We compute hyperdeterminants of hypermatrices whose indices
 belongs in a meet-semilattice and whose entries depend only of the greatest lower bound of the indices. One shows that an
 elementary expansion of such a polynomial allows to generalize a
 theorem of Lindstr\"om to higher-dimensional determinants. And we gave as an application
  generalizations of some results due to
Lehmer, Li and Haukkanen.
\end{abstract}
{\tiny {\it keywords:}Hyperdeterminants, meet semilattices,
GCD-matrices,
multiplicative functions. \\
{\it AMS:} 15A15, 15A69, 06A12.}
\section{Introduction}
Since the end of the nineteen century, it is known that some
determinants, with entries depending only of the $\gcd$ of the
indices,  factorize. Readers interested in the story of the problem
 can refer to \cite{Krat2}  and \cite{AST}. In 1876,
  Smith \cite{Smith}  evaluated the determinant of a GCD matrix whose entries belong
to a factor closed set (\ie\ a set which contains all the factors of
its elements) as a product of Euler's totient. The interest of this
computation lies in its links with arithmetic functions \cite{Apo}
and in particular multiplicative functions (see \cite{lasc,notelasc}
for interesting remarks about the last notion). During the last
century, many generalizations of Smith's theorem have been
investigated. One of the ways to extend his result consists in
changing the set of the indices of the matrices. In 1990, Li
\cite{Li}  gave the value of GCD determinant for an arbitrary set of
indices.
 Beslin and Ligh \cite{BL2} shown that such a determinant factorizes
 when the indices belongs to a $\gcd$-closed set
(\ie\ a set which contains the $\gcd$ of any pairs of its elements)
as a product of certain functions evaluated in terms of Euler's
totient. The fact that these determinants factorize can be seen as a
corollary of a very elegant theorem due to Lindstr\"om \cite{Lin1}
which evaluated the determinant of the GCD-matrix whose indices are
chosen in a meet semilattice, \ie\ a poset such that each pair
admits a greatest lower bound. Another way to generalize Smith's
result consists in computing multidimensional analogous. In 1930
Lehmer gave \cite{Leh} the first multi-indexed version of  Smith's
determinant. Other related computation are collected in
\cite{So1,So2}. More recently, Haukkanen \cite{Hau}  generalized the
results of  Beslin and Ligh  \cite{BL2} and  Li \cite{Li}) to
hyperdeterminants.

We will see in Section \ref{Hyperdet} that the main trick for
computing these multidimensional determinants consists in expanding
it as a sum of (classical $2$-way) determinants. In the aim to
highlight this method, we apply it to a more general object
$\Det_{\cal F}$.\\ In Section \ref{Meet}, we recall shortly a
classical technic and give a slight generalization of Lindstr\"om's
Theorem. As a consequence, we give a multidimensional analogue of
Lindstr\"om's Theorem.\\ In Section \ref{Minors}, one investigates
minors of meet hypermatrices and generalizes two theorems due to
Haukkanen \cite{Hau}.

\section{Hyperdeterminants and $\cal F$-determinants}\label{Hyperdet}

The question of extending the notion of determinant to higher
dimensional arrays has been raised by Cayley \cite{Ca1,Ca2} few
after he introduced the modern notation  as square arrays
\cite{Ca0}. The simplest generalization is defined for a $k$th order
tensor on an $n$-dimensional space $M=(M_{i_1,\cdots,i_k})_{1\leq
i_1,\cdots,i_k\leq n}$ by the alternated sum
\[\Det M=\frac1{n!}\sum_{\sigma=(\sigma_1,\cdots,\sigma_k)\in\S_n^k}\sign(\sigma)
M^\sigma,\] where
$\sign(\sigma)=\sign(\sigma_1)\cdots\sign(\sigma_k)$,
$M^\sigma=M_{\sigma_1(1)\dots\sigma_k(1)}\cdots
M_{\sigma_1(n)\dots\sigma_k(n)}$ and $\S_n$ is the symmetric group.
A straightforward computation
gives $\Det M=0$ if $k$ is odd.\\
 For any $k$ (even if $k$ is odd), one  defines
 the polynomial
\[ \Det_1 M=\sum_{\sigma=(\Id,\sigma_2,\cdots,\sigma_k)\in\S_n^k}
\sign(\sigma) M^\sigma.
\]
When $k$ is even the two notions coincide but for $k$ odd, only
$\Det_1$ does not vanish. This is a special case of the
"less-than-full-sign" determinant theory due to Rice \cite{Rice}.\\
  Let us denote by
${\cal F}$  a map from $\S_n^{k-2}$ to a commutative ring. One
defines a more general object, which will be called $\cal F$-{\it
determinant  of $M$}  by
\[\Det_{\cal F}(M)=\sum_{\sigma=(\sigma_2,\cdots,\sigma_{k})\in\S_n^{k-1}}
\sign(\sigma_2){\cal F}(\sigma_3,\dots,\sigma_k)
\prod_iM_{i\sigma_2(i)\dots\sigma_k(i)}.\] It exists an elementary
identity  which  consists in expanding the
 $\cal F$-determinant as a sum of $(n!)^{k-2}$ classical ($2$-way)
 determinants.
 \begin{lemma}{(Determinantal expansion})\\\label{splitting}
One has
\[
\Det_{\cal F} M=\sum_{\sigma_3,\dots,\sigma_k}{\cal
F}(\sigma_3,\dots,\sigma_k)\det(M^{\sigma_3,\dots,\sigma_k}),\]
where $M^{\sigma_3,\dots,\sigma_k}$ denotes the $n\times n$ matrix
such that
$M^{\sigma_3,\dots,\sigma_k}_{i,j}=M_{i,j,\sigma_3(i),\dots,\sigma_k(i)}$.
\end{lemma}
\Proof It suffices to remark that \[ \Det_{\cal
F}M=\sum_{\sigma_3,\dots,\sigma_k}{\cal
F}(\sigma_3,\dots,\sigma_k)\sum_{\sigma_2}\sign(\sigma_2)\prod_iM_{i\sigma_2(i)\sigma_3(i)\cdots\sigma_k(i)}
\]
\qed
 One of the most important property of  hyperdeterminants is the invariance under
the action of $k$ copies of the special linear group. It is a very
classical result which can be recover as a straightforward
consequence of the following proposition.
\begin{proposition}
The polynomial $\Det_{\cal F}M$ is invariant under the action of
linear group on $M$ in the following sense
\begin{equation}\label{Det_Finv}
\Det_{\cal F}g.M=\det g\Det_{\cal F}M,
\end{equation}
where
\[g.M=\left(\sum_{1\leq j_2 \leq n}g_{i_2j_2}M_{i_1,j_2,i_3\cdots,i_k}\right)_{1\leq i_1,\cdots,i_k\leq
n}.\]
\end{proposition}
\Proof By applying Lemma \ref{splitting} to $g.M$, one gets
\[\begin{array}{rl}
\Det_{\cal F} g.M=&\displaystyle\sum_{\sigma_3,\dots,\sigma_k}{\cal
F}(\sigma_3,\dots,\sigma_k)\det(g.M^{\sigma_3,\dots,\sigma_k})\\
=&\displaystyle\sum_{\sigma_3,\dots,\sigma_k}{\cal
F}(\sigma_3,\dots,\sigma_k)\det g
\det(M^{\sigma_3,\dots,\sigma_k})\\
=&\det g\Det_{\cal F}M.\end{array}\] \qed
\section{Hyperdeterminants of meet hypermatrices}\label{Meet}
\subsection{Meet semilattice}

Consider a partially ordered finite set $L$ so that every pairs
$(x,y)\in L^2$ has a greatest lower bound denoted by $x\wedge y$.
 Such a poset is called a meet semilattice. One defines classically
its $\zeta$ function by \[\zeta(x,y)=\left\{\begin{array}{ll}
1&\mbox{ if } x\leq y,\\ 0 &\mbox{ otherwise}.
\end{array}\right.\]
Its M\"obius function is the inverse of the zeta function and can be
computed by the induction
\[
\mu(x,y)=\left\{\begin{array}{ll} 1&\mbox{ if } x=y,\\
-\sum_{x\leq z<y}\mu(z,y)&\mbox{ if }x<y,\\
0&\mbox{ in the other cases}.
\end{array}\right.
\]
 If $F$ and $f$ verify the equality
\begin{equation}
\label{mo1} F(x)=\sum_{y\leq x}f(y)=\sum_{y\in
L}\zeta(y,x)f(y),
\end{equation}
then, one has
\begin{equation}\label{mobiusinv}
f(x)=\sum_{y\in L}\mu(y,x)F(y)=(F\star \mu)(x)
\end{equation}
where the symbol $\star$ means the convolution product.
\subsection{Lindstr\"om Theorem} The factorization properties of the GCD determinants
are the consequence of the semilattice structure of the integers
with respect to divisibility and can be stated in a more general
way. The manipulations of the identities (\ref{mo1}) and
(\ref{mobiusinv}) are the keys of the proof of Lindstr\"om's Theorem
\cite{Lin1}. We recall its  proof in a very slightly more general
version.

For each $x\in L$, one considers a fixed element $z_x\leq x$. Let
$F_x$ be a function from $L$ to $\C$ (or more generally to a
commutative ring). Let $M$ be the matrix defined by \[
M=\left(F_x(z_x\wedge y)\right)_{x,y\in L}.
\]
Remark that it suffices to define $F_x(z)$ only when $z\leq x$. In
particular, one can suppose that $F_x(z)=F(z,x)$ is an incidence
function, \ie\ $F(x,z)=0$ unless $x\leq z$. One has
\[
\begin{array}{ll}F_x(z_x\wedge y)=&\displaystyle\sum_{z\in
L}\zeta(z,z_x)\zeta(z,y)f_x(z),
\end{array}\]
where $f_x(z)=\sum_{y\in L}\mu(y,z)F_x(y)$. Hence, $\det M$
factorizes as the product
\[\ \det M= \det \Phi.\det Z, \] where
$\Phi=\left(\zeta(y,z_x)f_x(y)\right)_{x,y\in L}$ and
$Z=\left(\zeta(x,y)\right)_{x,y\in L}$. As $\Phi$ and $Z$ are
triangular, $\det Z=1$ and
\begin{equation}\label{cesarogen}
\det\Phi=\prod_x\zeta(x,z_x)f_x(x)=\left\{\begin{array}{ll}
 \displaystyle\prod_xf_x(y)&\mbox{ if }z_x=x \mbox{ for each }x,\\
 0&\mbox{otherwise}.
\end{array}\right.
\end{equation}
Then one obtains Lindstr\"om's Theorem.
\begin{theorem}(Lindstr\"om \cite{Lin1})\label{meetdet}
\begin{equation}
\det\left(F_x(z_x\wedge y)\right)_{x,y\in
L}=\left\{\begin{array}{ll}
 \prod_xf_x(x)&\mbox{ if }z_x=x \mbox{ for each }x,\\
 0&\mbox{otherwise}.
\end{array}\right.
\end{equation}
\end{theorem}
Note that, the original Lindst\"om 's Theorem deals with the case
where $z_x=x$ for each $x$. Furthermore, equality (\ref{cesarogen})
generalizes a lemma of Cesaro.
\begin{lemma}\label{Cesaro}(Cesaro)\\Denote by $\gcd_{m}(n)=\gcd(m,n)$.
One has,\begin{equation}(\mu\ast(f\circ{\rm
gcd}_m))(n)=\left\{\begin{array}{ll}(f\ast \mu)(n)&\mbox{ if
}m=n,\\0&\mbox{ otherwise},\end{array}\right.\end{equation} where
$\ast$ is the Dirichlet convolution and $\circ$ is the composition
of functions.
\end{lemma}
\subsection{Linstr\"om's Theorem for $\cal F$-determinants}
Lindst\"om's Theorem  can be extended to $\cal F$-determinants.
\begin{theorem}\label{Linhypdet} (Lindstr\"om's theorem for $\cal F$-determinants)\\
If $L=\{x_1,\dots, x_n\}$ denotes a meet semilattice, one has
\[
\Det_{\cal F}\left(F_{x_{i_1}}(z_{x_{i_1}}\wedge\cdots\wedge
x_{i_k})\right)=\left\{\begin{array}{ll}{\cal
F}(\Id,\cdots,\Id)\displaystyle\prod_xf_x(x)&\mbox{if }z_x=x\mbox{
for each } x\\0&\mbox{otherwise}\end{array}\right.\]
\end{theorem}
\Proof Lemma \ref{splitting} gives
\[\begin{array}{l}
\Det_{\cal F}\left(F_{x_1}(z_{x_1}\wedge\cdots\wedge
x_k)\right)_{x_1,\cdots,x_k\in L}=\\
 \displaystyle\sum_{\sigma_3,\dots,\sigma_k}{\cal
 F}(\sigma_3,\dots,\sigma_k)\det\left(F_{x_i}(z_{x_i}\wedge x_j\wedge x_{\sigma_3(i)}
\wedge \cdots\wedge x_{\sigma_k(i)}\right)_{i,j}\end{array}
\]
%and the result follows from  Lindst\"om's Theorem (Theorem
%\ref{meetdet}).
From Linstr\"om's Theorem (Theorem \ref{meetdet}), one has
\[\det\left(F_{x_i}(z_{x_i}\wedge x_j\wedge x_{\sigma_3(i)} \wedge
\cdots\wedge x_{\sigma_k(i)}\right)\neq 0\] if and only if for each
$x_i$ one has
\[z_{x_i}\wedge x_{\sigma_3(i)}\wedge\cdots\wedge x_{\sigma_k}(i)=x_i.\]
Equivalently,
\[z_{x_i}\wedge x_{\sigma_3(i)}\wedge\cdots\wedge x_{\sigma_k}(i)\geq x_i,\]
for eacu $i$. Hence,  $\sigma_3=\cdots=\sigma_k=Id$ and
$z_{x_i}=x_i$ for each $i$. The result follows.
 \qed
 \begin{example}
Consider the semilattice constituted with two elements $2\geq 1$.
The expansion of the $\goth F$-determinant gives
\[\small
\begin{array}{rcl}
\DetF(F_{i}(i\wedge j\wedge k\wedge l))&=&{\goth
F}(12,12)\left|\begin{array}{cc}F_{1}(1)&F_{1}(1)\\
F_{2}(1)&F_{2}(2)\end{array}\right|\\&&+\left({\goth
F}(12,21)+{\goth F}(21,12)+{\goth F}(21,21)\right)\left|\begin{array}{cc}F_{1}(1)&F_{1}(1)\\
F_{2}(1)&F_{2}(1)\end{array}\right|\\
&=&{\goth F}(12,12)F_{1}(1)(F_2(2)-F_2(1))\\
&=&{\goth F}(12,12)f_1(2)f_2(2)
\end{array}
\]
 \end{example}
\section{Minors of meet Hypermatrices}\label{Minors}
\subsection{Meet closed subsets}
Consider a meet closed subset $S$ of $L$ (\ie\, a subset closed
under the operation $\wedge$) and fix a linear extension $l=y_1\dots
y_n$ of $S$. Following the notations of \cite{AST}, we denote by
$x\trianglelefteq y_i$ the relation $x\leq y_i$ and $x\not\leq y_j$
for each $j<i$. Consider a pair of functions $f$ and $F$ verifying
\begin{equation}\label{eqshat}F(y_i)=\sum_{x\leq y_i\atop x\in L}f(x)\end{equation} and set $\hat
f(y_i)=\sum_{x\trianglelefteq y_i}f(x)$. One has the following
lemma.
\begin{lemma}\label{hat}
\begin{equation}\label{eqhat} F(y_i)=\sum_{y_k\leq y_i}\hat f(y_k).\end{equation}
\end{lemma}
\Proof Remarking that for each $x\in L$, it exists $i$ such that
$x\trianglelefteq y_i$ (it suffices to set $y_i=\min\{j|x\leq y_j\}$
and that $x\trianglelefteq y_i$ and $x\trianglelefteq y_j$ implies
$i=j$, we have
\[F(y_i)=\sum_{x\leq y_i\atop x\in L}f(x)=
\sum_{y_k\leq y_i}\sum_{x\trianglelefteq y_k}f(x)=\sum_{y_k\leq y_i}\hat f(y_k).\]
%This identity appears in \cite{AST} (Theorem 4.1, p 7)
%and its proof is based on a one-to-one correspondance between the
%terms in (\ref{eqhat}) and (\ref{eqshat}). We do not repeat it here.
\qed Note that this identity appears in \cite{AST} (Theorem 4.1, p
7). Hence, using Theorem \ref{Linhypdet} and Lemma \ref{hat}, one
generalizes a result by Altinisik, Sagan and Tuglu (\cite{AST},
Theorem 4.1 p 7).
\begin{corollary}\label{Hauklattice}
\begin{equation}\nonumber
\Det_{\cal F}(F_{y_{i_1}}(y_{i_1}\wedge\dots\wedge y_{i_k}))={\cal
F}(\Id,\cdots,\Id)\prod_{i=1}^n\left(\sum_{x_1\trianglelefteq
y_i}\sum_{x_2\in L}\mu(x_1,x_2)F_{y_i}(x_2)\right).
\end{equation}
\end{corollary}
\begin{example}
Consider the semilattice $L$ given by its Hasse diagram
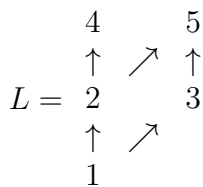
\[
L=\begin{array}{ccc}
4&&5\\\uparrow&\nearrow&\uparrow\\2&&3\\\uparrow&\nearrow&\\1&&\end{array}
\]
where $i\rightarrow j$ means $i\leq j$.
 The sublattice generated by $2$, $4$ and $5$,
\[
S=\begin{array}{ccc} 4&&5\\\uparrow&\nearrow&\\2&&\end{array}
\]
 is meet closed  and
\[
\begin{array}{l}\hat f_2(2)=f_2(2)+f_2(1)\\\hat f_4(4)=f_4(4)\\\hat f_5(5)=f_5(5)+f_5(3). \end{array}
\]
Hence,
\[\small
\begin{array}{rcl}
\DetF\left(F_i(i\wedge j\wedge k)\right)_{i,j,k\in S}&=& {\goth
F}(123)\left|\begin{array}{ccc}F_2(2)&F_2(2)&F_2(2)\\
F_4(2)&F_4(4)&F_4(2)\\F_5(2)&F_5(2)&F_5(5)\end{array}\right|\\&&+
 {\goth
F}(213)\left|\begin{array}{ccc}F_2(2)&F_2(2)&F_2(2)\\
F_4(2)&F_4(2)&F_4(2)\\F_5(2)&F_5(2)&F_5(5)\end{array}\right|\\&&+
{\goth
F}(321)\left|\begin{array}{ccc}F_2(2)&F_2(2)&F_2(2)\\
F_4(2)&F_4(2)&F_4(2)\\F_5(2)&F_5(2)&F_5(2)\end{array}\right|\\&&+
({\goth F}(132)+{\goth F}(231)+{\goth
F}(312))\left|\begin{array}{ccc}F_2(2)&F_2(2)&F_2(2)\\
F_4(2)&F_4(2)&F_4(2)\\F_5(2)&F_5(2)&F_4(2)\end{array}\right|\end{array}\]

The permutation $123$ is the only one having a non zero contribution
in this sum. Hence,
\[\begin{array}
{rcl} \DetF\left(F_i(i\wedge j\wedge k)\right)_{i,j,k\in S}
&=&{\goth F}(123)F_2(2)(F_4(4)-F_4(2))(F_5(5)-F_5(2))\\&=&
 {\goth F}(123)(f_2(2)+f_2(1))f_4(4)(f_5(5)+f_5(3))\\
 &=&{\goth F}(123)\hat f_2(2)\hat f_4(4)\hat f_5(5)
\end{array}
\]
\end{example}
\begin{remark}
If $L$ is the semilattice structure of the integers with respect to
divisibility. By the specialization ${\cal
F}(\sigma_3,\dots,\sigma_k)=\sign(\sigma_3)\dots\sign(\sigma_k)$,
one recovers the computations of Lehmer \cite{Leh} as a special case
of Theorem \ref{Linhypdet} and the result of Haukkanen (\cite{Hau}
Theorem 1. p 56) from Corollary \ref{Hauklattice}.
\end{remark}
\subsection{Factor closed subsets} Let $S$ be a factor closed subset
of $L$. Then, $f=\hat f$ and

\begin{equation}\label{Lehm}
\Det_{\cal F}\left(F_{x_1}(x_1\wedge\cdots\wedge
x_k)\right)_{x_1,\dots,x_k\in S}={\cal
F}(\Id,\cdots,\Id)\prod_xf_x(x).\end{equation}

As special cases of equality (\ref{Lehm}), one recovers  Lehmer's
identities \cite{Leh} and the original result of Smith \cite{Smith}.
\subsection{General case} Let
$X=\{x_1,\dots,x_n\}$ be a subposet of a meet semi-lattice $L$. We
will denote by $\overline X=\{x_1,\dots,x_n,x_{n+1},\dots,x_m\}$ the
smallest factor-closed subset of $L$ containing $X$.
 The aim of this section consists in investigating the $F$-determinant
\begin{equation}DF(X):=\DetF\left(F_{x_{i_1}}(z_{x_{i_1}}\wedge
x_{i_2}\wedge\cdots\wedge x_{i_k})\right)_{1\leq i_1,\dots,i_k\leq
n}\end{equation} where for each $x\in X$, $z_x$ denotes a fixed
element of $\overline X$ such that $z\leq x$. As in the previous
section, the result follows from the case $k=2$ and Proposition
\ref{splitting}.\\
Let us consider first the determinant $\det(F_{x_i}(z_{x_i}\wedge
x_j))_{1\leq i,j\leq n}$. The functions $F_{x_i}$ can be chosen such
that $F(x,y)$ is an incidence function. The set $\overline X$ being
closed by factors, the functions $f_x$ and $\hat f_x$ are equal.
Hence,
\begin{equation} F_{x_i}(z_{x_i}\wedge
x_j)=\sum_{k=1}^{n+m}C_{x_i,x_k}\zeta(x_k,x_j)
\end{equation}
where \begin{equation}
C_{x,y}=f_x(y)\zeta(y,z_x)=\left\{\begin{array}{ll}f_x(y)&\mbox{if }
y\leq z_x\\0&\mbox{otherwise}.\end{array}\right.\end{equation} One
has
\begin{proposition}\label{Li}
\begin{equation}\det(F_{x_i}(x_i\wedge x_j)) =\sum_{1\leq k_1<\dots<k_n\leq
n+m}\det\left(C_{x_i,x_{k_j}}\right)_{1\leq i,j\leq
n}\det\left(\zeta(x_{k_i},x_{j})\right)_{1\leq i,j\leq n}.
\end{equation}
\end{proposition}
\Proof By multi-linearity, one obtains \[ \small\begin{array}{ll}
\det(F_{x_i}(x_i\wedge
x_j))&=\det\left(\displaystyle\sum_{k=1}^{n+m}C_{x_i,x_k}\zeta(x_k,x_j)\right)
\\&=\displaystyle\sum_{1\leq k_1,\dots,k_n\leq
n+m}\det\left(C_{x_i,x_{k_j}}\zeta(x_{k_j},x_j)\right)\\&=
\displaystyle\sum_{1\leq k_1,\dots,k_n\leq
n+m}\det\left(C_{x_i,x_{k_j}}\right)\prod_i\zeta(x_{k_i},x_i)\\
&=\displaystyle\sum_{1\leq k_1<\dots<k_n\leq
n+m}\det\left(C_{x_i,x_{k_j}}\right)_{1\leq i,j\leq
n}\det\left(\zeta(x_{k_i},x_{j})\right)_{1\leq i,j\leq
n}\end{array}\] \qed
%\begin{example}
%Consider the semilattice
%\[
%L:=\begin{array}{cccccc}c&&&d&&\\&\nwarrow&\nearrow&&\\&&b&&&e\\&&&\nwarrow&\nearrow&\\&&&a&&
%\end{array}
%\]
%And the sub-poset $X$ composed by the elements $c$, $d$ and $e$, and
%suppose $z_c=b$, $z_d=d$ and $z_e=e$. We want to compute the
%determinant
%\[
%\det(F_x(z_x\wedge y\wedge t))_{x,y,t\in
%X}=\left|\begin{array}{ccc}f_c(b)+f_c(a)&f_c(b)+f_c(a)&f_c(a)\\
%f_d(b)+f_d(a)&f_d(d)+f_d(b)+f_d(a)&f_d(a)\\f_e(a)&f_e(a)&f_e(e)+f_e(a)
%\end{array}\right|
%\]
%if one expands this determinant by multilinearity and one collects ,
¨%one gets
\begin{example}
Let us consider the semilattice $L$ whose Hasse diagram is
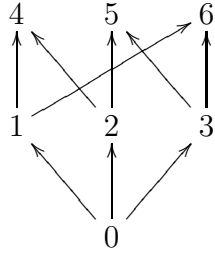
\[
\newdir{ >}{{}*!/-4mm/@{>}}
\xymatrix@R=10mm@L=10mm{
  4 & 5 &  6 \\
 1 \ar@{->}[u] \ar@{->}[urr]   &2 \ar@{->}[u] \ar@{->}[ul] & 3 \ar@{->}[ul] \ar@{->}[u]
 \\
 &  0 \ar@{->}[ur] \ar@{->}[ul] \ar@{->}[u]&}
%& \SMQSym \ar@{<-<}[dl] \ar@{->>}[ddr]\ar@{<-<}[rrr]&&&\MQSym\ar@{<-<}[dl] \ar@{->>}[ddr]\\
%\SMRSym\ar@{->>}[ddr]\ar@{<-<}[rrr]& &&\MRSym\ar@{->>}[ddr]\\
%&&\SMCSym\ar@{<-<}[dl]\ar@{<-<}[rrr]&&&\MCSym\ar@{<-<}[dl]\\
%  &\SMSym\ar@{<-<}[rrr]&&&\MSym}
\]
and $X=\{4,5,6\}$. We set $z_4=1$, $z_5=2$ and $z_6=6$ . Consider
the determinant
\[\small\begin{array}{rcl}
\det(F_i(z_i\wedge j\wedge k))_{i,j,k\in X}&=&\left|\begin
{array}{ccc} f_4(1)+f_4(0)&f_4(0)&f_4(1) + f_4(0)\\f_5(2)+f_5(0)&
f_5(2)+f_5(0)&f_5(0)
\\ f_6(1)+ f_6(0)&f_6(3)+
f_6(0)&f_6(6)+ f_6(3)+ f_6(1)+f_6(0)
\end {array}\right|\\&=&
-f_4(0)f_6(1){\it f_5}(2)+{\it f_4}(1){\it f_6}_{{0 }}{\it
f_5}(2)+{\it f_4}(1){\it f_5}(0){\it f_6}(3)\\&&+{\it f_4}_{
{1}}{\it f_5}(0){\it f_6}(6)+{\it f_4}(0){\it f_5}(2){\it f_6}_
{{3}}+2\,{\it f_4}(1){\it f_5}(2){\it f_6}(3)\\&&+{\it f_4}(1){ \it
f_5}(2){\it f_6}(6)
\end{array}
\]
Using the multilinearity of $\det$, one recovers the expression
given by Proposition \ref{Li}
\[\small\begin{array}{l}
\det(F_i(z_i\wedge j))=\left|\begin {array}{ccc} f_4(0)&f_4(1)&0\\
f_5(0)&0& f_5(2)
\\f_6(0)&f_6(1)&0\end {array}\right|\cdot\left|
\begin {array}{ccc} 1&1&0\\ 1&0&1\\
1&1&0\end {array}\right|+
 \left|\begin {array}{ccc} {\it f_4}(0)&
{\it f_4(1)}&0\\ {\it f_5}(0)&0&0\\ f_6(0)&f_6(1)&{\it
f_6}(3)\end {array}\right|\cdot\left|\begin {array}{ccc} 1&1&0\\
1&0&0
\\ 1&1&1\end {array}\right|\\
+\left|\begin {array} {ccc} {\it f_4}(0)&{\it f_4}(1)&0\\
f_5(0)&0&0\\ f_6(0)&f_6(1)&{\it f_6}(6)
\end {array}\right|\cdot\left|\begin {array}{ccc} 1&1&0
\\ 1&0&0\\ 1&1&1\end {array}\right|+
\left|\begin {array}{ccc} {\it f_4}(0)&0&0\\
f_5(0)&{\it f_5}(2)&0\\ f_6(0)&0&{ \it f_6}(3)\end
{array}\right|\cdot\left|\begin {array}{ccc} 1&0&0
\\ 1&1&0\\ 1&0&1\end {array}\right|\\
+\  \left|\begin {array}{ccc} {\it f_4}(1)&0&0
\\ 0&{\it f_5}(2)&0\\ f_6(1)&0&{\it f_6}
(3)\end {array}\right|\cdot\left|\begin {array}{ccc} 1&0&0
\\ 0&1&0\\ 1&0&1\end {array}\right|+
\left|\begin {array}{ccc} {\it f_4}(1)&0&0
\\ 0&{\it f_5}(2)&0\\ f_6(1)&0&{\it f_6}
(6)\end {array}\right|\cdot\left|\begin {array}{ccc} 1&0&0
\\ 0&1&0\\ 1&0&1\end {array}\right|\end{array}
\]

\end{example}

More generally, one has a multi-indexed version of Proposition
\ref{Li}.
\begin{theorem}\label{LiGen}
\begin{equation}\label{DF}
\begin{array}{r} DF(X)=\displaystyle\sum_{1\leq k_1<\cdots<k_n\leq
n+m}\DetF\left(f_{x_{i_1}}(x_{k_{i_2}})\zeta(x_{k_{i_2}},z_{x_{i_1}}\wedge
x_{i_3}\wedge\cdots\wedge
x_{i_k})\right)\times\\\times\det\left(\zeta(x_{k_i},x_j)\right)
\end{array}\end{equation}
\end{theorem}
\Proof We use Lemma \ref{splitting} to expand $DF(X)$ and we obtain
\[DF(X)=\sum_{\sigma_1,\dots,\sigma_k}{\goth F}(\sigma_1,\dots,\sigma_k)
\det\left(F_{x_i}(z_{x_i}\wedge x_j\wedge
x_{\sigma_3(i)}\wedge\cdots\wedge x_{\sigma_k(i)})\right).\] Now, by
Proposition \ref{Li}, one gets
\[\begin{array}{rcl}
DF(X)&=&\displaystyle    \sum_{\sigma_1,\dots,\sigma_k}{\goth
F}(\sigma_1,\dots,\sigma_k)\displaystyle\sum_{1\leq
k_1<\cdots<k_n\leq
n+m}\det\left(\zeta(x_{k_i},x_j)\right)\times\\&&\times
\det\left(f_{x_i}(x_{k_j}\zeta(x_{k_j},z_{x_i}\wedge
x_{\sigma_3(i)}\wedge\cdots\wedge x_{\sigma_k(i)}))\right)\\
&=&\displaystyle\sum_{1\leq k_1<\cdots<k_n\leq
n+m}\det\left(\zeta(x_{k_i},x_j)\right)\displaystyle
\sum_{\sigma_1,\dots,\sigma_k}{\goth
F}(\sigma_1,\dots,\sigma_k)\times\\&&\times
\det\left(f_{x_i}(x_{k_j}\zeta(x_{k_j},z_{x_i}\wedge
x_{\sigma_3(i)}\wedge\cdots\wedge x_{\sigma_k(i)}))\right)\\
&=&\displaystyle\sum_{1\leq k_1<\cdots<k_n\leq
n+m}\DetF\left(f_{x_{i_1}}(x_{k_{i_2}})\zeta(x_{k_{i_2}},z_{x_{i_1}}\wedge
x_{i_3}\wedge\cdots\wedge
x_{i_k})\right)\times\\&&\times\det\left(\zeta(x_{k_i},x_j)\right)
\end{array}\]
This ends the proof. \qed Let set $F_{x_1}=F_{x_2}=\cdots=F_{x_n}=F$
then for $1\leq k_1,k_2,\dots,k_n\leq n$, we get
% \[\begin{array}{ll}
%\det\left(C_{x_i,x_{k_j}}\right)&=\det\left(f(x_{k_j})\zeta(x_{k_j},x_i)\right)\\&=
%\det\left(\zeta(x_{k_i},x_j)\right)\displaystyle\prod_if(x_{k_i}).\end{array}
%\]
% Again, by Lemma \ref{splitting}, one obtains the following
%corollary.
\begin{corollary}\label{GenHauk}\[
\begin{array}{r} \DetF\left(F(z_{x_{i_1}}\wedge x_{i_2}\wedge\cdots\wedge
x_{i_k})\right)=\displaystyle\sum_{1\leq k_1<\cdots<k_n\leq
n+m}\prod_if(x_{k_i})\times\\\times
\DetF\left(\zeta(x_{k_{i_2}},z_{x_{i_1}}\wedge
x_{i_3}\wedge\cdots\wedge
x_{i_k})\right)\det\left(\zeta(x_{k_i},x_j)\right)
\end{array}\]
\end{corollary}
\Proof By applying the equality,
 \[
\small
 \begin{array}{l}
 \DetF\left(f(x_{k_{i_1}})\zeta(x_{k_{i_2}},z_{x_{i_1}}\wedge
 x_{i_3}\wedge\dots\wedge x_{i_k})\right)=\\
 \displaystyle\sum_{\sigma_3,\dots,\sigma_k}{\goth
 F}(\sigma_3,\cdots,\sigma_k)\det\left(f(x_{k_{i}})\zeta(x_{k_{j}},z_{x_{i}}\wedge
 x_{\sigma_1(i)}\wedge\dots\wedge x_{\sigma_k(i)})\right)=\\ \displaystyle\prod_if(x_{k_i})
\DetF\left(\zeta(x_{k_{i_2}},z_{x_{i_1}}\wedge
x_{i_3}\wedge\cdots\wedge x_{i_k})\right).
 \end{array}
 \]
 to identity \ref{DF}, one obtains the result.
\qed
\begin{remark} Assume that $L$ is the integer lattice. Then, if $z_x=x$ for
each $x\in X$ and ${ \goth
F}(\sigma_3,\cdots,\sigma_k)=\sign(\sigma_3)\dots\sign(\sigma_k)$ in
Corollary \ref{GenHauk}, one recovers Theorem 2 in \cite{Hau}.
Moreover, Proposition \ref{Li} generalizes the  result of
Li\cite{Li}.\end{remark}

\end{document}